\documentclass[a4paper,12pt]{amsart}
\newcommand{\Qbar}{\overline{\mathbb Q}}
\newcommand{\ZZ}{{\mathbb Z}}
\newcommand{\QQ}{{\mathbb Q}}
\newcommand{\cO}{{\mathcal O}}
\newcommand{\cE}{{\mathcal E}}
\newcommand{\cT}{{\mathcal T}}
\newcommand{\cC}{{\mathcal C}}

\usepackage{amssymb}

\newtheorem{theorem}{Theorem}

\newtheorem{e-proposition}[theorem]{Proposition}

\newtheorem{e-definition}[theorem]{Definition\rm}
\newtheorem{remark}{\it Remark\/}

\title{On the ubiquity of trivial torsion on elliptic curves}
\author[Enrique Gonz\'alez-Jim\'enez]{Enrique Gonz\'alez-Jim\'enez}
\address{Universidad Aut{\'o}noma de Madrid, Departamento de Matem{\'a}ticas and Instituto de Ciencias Matem{\'a}ticas (CSIC-UAM-UC3M-UCM), Madrid, Spain}
\email{enrique.gonzalez.jimenez@uam.es}
\author{Jos\'e M. Tornero}
\address{Departamento de \'Algebra, Universidad de Sevilla. P.O. 1160. 41080 Sevilla, Spain.}
\email{tornero@us.es}
\date{\today}
\begin{document}
\maketitle

\begin{abstract}
The purpose of this paper is to give a {\it down--to--earth} proof of the well--known fact that a randomly chosen elliptic curve over the rationals is most likely to have trivial torsion. \\
\end{abstract}
\section{Introduction}
Let us consider an elliptic curve $E$, defined over the rationals
and written in short Weierstrass form
\begin{equation}\label{eq_short}
E: Y^2 = X^3 + AX + B, \quad A,B \in \ZZ.
\end{equation}
We will use the standard notations for:

\begin{itemize}
\item $\Delta = -16(4A^3+27B^2) \neq 0$, the discriminant of $E$;
\item $E(\QQ)$, the finitely generated abelian group of rational points on $E$, and 
\item $\cO$, the identity element of $E(\QQ)$. 
\end{itemize}

Given $P \in E(\QQ)$, we will also write as customary $[m]P$ for the point resulting after adding $m$ times  $P$.

The problem of computing the torsion of $E(\QQ)$ has been solved in a
lot of very efficient ways \cite{Cremona,Doud,GOT}, and most
computer packages (say \verb+Maple-Apecs+, \verb+PARI/GP+, \verb+Magma+ or \verb+Sage+)
calculate the torsion of curves with huge coefficients in very
few seconds. The major result which made this possible (along with
others, like the Nagell--Lutz Theorem (\cite{Nagell},\cite{Lutz}) or the embedding theorem for good
reduction primes (see, for example, \cite[VIII.7]{Silverman} or  \cite[Chap. 5]{Husemoller})) was Mazur's Theorem
\cite{MazurIHES,Mazur} who listed the fifteen possible torsion
groups. 

In the above papers, it is proved that the possible structures of the torsion group of $E(\QQ)$ are
$$ 
\ZZ / n \ZZ \mbox{ for } n=2,\dots,10,12, \quad \mbox{ or } \quad \ZZ / 2 \ZZ\times \ZZ / 2n \ZZ \mbox{ for } n=1,\dots,4.
$$ 

Besides, the fifteen of them actually happen as torsion subgroups of elliptic curves. Notice that thanks to the above theorem, the possible prime orders for a torsion point defined over $\QQ$ are $2,3,5$ or $7$.

Let $p$ be a prime number and let $E[p]$ be the group of points of order $p$ on $E(\Qbar)$, where $\Qbar$ denotes an algebraic closure of $\QQ$. The action of the absolute Galois group ${\rm G}_\QQ={\rm Gal}(\Qbar/\QQ)$ on $E[p]$ defines a mod $p$ Galois representation
$$
\rho_{E,p}:{\rm G}_\QQ\rightarrow {\rm Aut}(E[p])\cong {\rm GL}_2(\mathbb{F}_p).
$$
Let $\QQ(E[p])$ be the number field generated by the coordinates of the points of $E[p]$. Therefore, the Galois extension $\QQ(E[p])/\QQ$ has Galois group
$$
{\rm Gal}(\QQ(E[p])/\QQ)\cong \rho_{E,p}({\rm G}_\QQ).
$$
The prime $p$ is called exceptional for $E$ if $\rho_{E,p}$ is not surjective. If $E$ has complex multiplication then any odd prime number is excepcional. On the other hand, if $E$ does not have complex multiplication then Serre \cite{Serre2} proved that $E$ has only finitely many exceptional primes. 

Duke  \cite{Duke} proved that {\em almost all} elliptic curves over $\QQ$ have no exceptional primes. More precisely, given an elliptic curve $E$ in a short Weierstrass form as in (\ref{eq_short}), the height of the elliptic curve is defined as 
$$
H(E) = \max(|A|^3,|B|^2).
$$ 

Let $M$ be a positive integer, and let $\mathcal{C}_H(M)$ be the set of elliptic curves $E$ with $H(E)\le M^6$. For any prime $p$ denote by $\cE_p(M)$ the set of elliptic curves $E\in \mathcal{C}_H(M)$ such that $p$ is an excepcional prime for $E$, and by $\cE(M)$ the union of $\cE_p(M)$ for all primes. Actually in both sets the elliptic curves were considered up to $\QQ$--isomorphisms. Duke then proved that
$$
\lim_{M \to \infty} \frac{|\cE(M)|}{|\mathcal{C}_H(M)|} = 0.
$$
His proof is based on a version of the Chebotarev density theorem, and uses a two-dimensional large sieve inequality together with results of Deuring, Hurwitz and Masser-W\"ustholz.

Duke also conjectured the following fact, later proved by Grant \cite{Grant}
$$
 \left|\cE(M) \right| \sim c \sqrt{M}.
$$
Being a bit more precise, Grant showed that, in order to efficiently estimate $\left| 
\cE(M) \right|$, only $\cE_2(M)$ and $\cE_3(M)$ had to be actually taken 
into account.

\

Now recall that there is a tight relationship between exceptional primes and torsion orders, because if there is a point of order $p$, then $p$ is an exceptional prime \cite{Serre2}. Our aim is then giving a down-to-earth proof of the fact that {\em almost all} elliptic curves over $\QQ$ have trivial torsion, motivated by Duke's paper. 

We will use in order to achieve this the characterization of torsion structures given in \cite{GT,GT2}, Mazur's Theorem \cite{MazurIHES,Mazur}; and a theorem by Schmidt \cite{Schmidt} on Thue inequalites. Note that we have used a different height notion, more naive in some sense, but nevertheless better suited for our purposes. 

Let us change a bit the notation and let us call 
$$
E_{(A,B)}: \; Y^2=X^3+AX+B
$$ 
and, provided $\Delta \neq 0$, we will denote by $E_{(A,B)}(\QQ)[m]$ the group of 
points $P \in E_{(A,B)}(\QQ)$ such that $[m]P = \cO$. Let us write as well
\begin{eqnarray*}
\cC(M) &=& \{ (A,B) \in \ZZ^2 \; | \; \Delta = -16(4A^3+27B^2) \neq 0, \; \; |A|,|B| \leq M \}.\\
\cT_p(M) &=& \{ (A,B) \in \cC(M) \; | \; E_{(A,B)}(\QQ)[p] \neq \{ \cO \} \}. \\
\cT (M) &=& \bigcup_{p \; \mbox{\scriptsize prime }} \cT_p (M)
\end{eqnarray*}

Our version of Duke's result is then as follows.

\

\begin{theorem}\label{teorema} With the notations above,
$$
\lim_{M \to \infty} \frac{|\cT(M)|}{|\cC(M)| } = 0.
$$ 
\end{theorem}

\

The proof will lead to extremely coarse bounds for $|\cT_p(M)|$ which will be proved unsatisfactory in view of experimental data, which we will display subsequently. 

\section{Proof of Theorem \ref{teorema}.}

Recall that the possible prime orders of a torsion point defined over $\QQ$ are $2,3,5$ or $7$.

We will make extensive use of the parametrizations of curves with 
a point of prescribed order given in \cite{GT,GT2,Kubert}. These 
results have recently been proved useful in showing new properties 
of the torsion subgroup (see, for instance \cite{BI,GGT,Ingram}).

First, note that, for a given $A$ with $|A|\leq M$ there are, at
most, two possible choices for $B$ such that $\Delta =0$ (and
hence, the corresponding curve $E_{(A,B)}$ is not an elliptic curve).
Therefore
$$
|\cC(M)| \geq (2M+1)^2 - 2(2M+1) = 4M^2-1.
$$

Let us recall from \cite{GT} that a curve $E_{(A,B)}$ with
a point of order $2$ must verify that there exist $z_1,z_2 \in
\ZZ$ such that
$$
A = z_1-z_2^2, \quad\quad\quad B = z_1z_2.
$$
Therefore $z_1|B$ and for a chosen $z_1$, both $z_2$ and $A$ are determined. Hence, there is at most one pair in $\cT_2(M)$ for every divisor of $B$. 

We need now an estimate for the average order of
the function $d(x)$, the number of positive divisors of $x$. The simplest estimation is, 
probably, the one that can be found in \cite{HW},
$$
d(1)+d(2)+...+d(x) \sim x \log(x).
$$ 
Therefore, as $M$ tends to infinity,
$$
\left| \cT_2(M) \right| \leq \sum_{x=1}^M 2 d(x) + \sum_{x=1}^M 2 d(x) + 2M,
$$
taking into account that we need to consider both positive and negative divisors, the cases where $x \in \{-M,...,-1\}$ and the $2M$ curves with $B=0$. Hence $\left| \cT_2(M) \right|\sim c_2 M\log(M)$, where we can, in fact, take $c_2 = 4$.

As for points of order $3$ we can find in \cite{GT} a similar
characterization (a bit more complicated this time) based on the
existence of $z_1,z_2 \in \ZZ$ such that 
$$
A = 27z_1^4+6z_1z_2, \quad\quad\quad B = z_2^2-27z_1^6.
$$
Analogously $z_1|A$ and, once we fix such a divisor, $z_2$ is necessarily given by
$$
z_2 = \frac{A-27z_1^4}{6z_1},
$$
which implies that again there is at most one pair in $\cT_3(M)$
for every divisor of $A$. Hence, as $M$ tends to infinity
$$
\left| \cT_3(M) \right| \leq c_3 M\log(M),
$$
and again $c_3=4$ suits us.

Points of order $5$ and $7$ need a similar, yet slightly different
argument. From \cite{GT2} we know that if there is a point of order $5$
in $E_{(A,B)}(\QQ)$, then there must exist $p,q \in
\ZZ$ verifying:
\begin{eqnarray*}
A &=& -27(q^4-12q^3p+14q^2p^2+12p^3q+p^4), \\
B &=& 54(p^2+q^2)(q^4-18q^3p+74q^2p^2+18p^3q+p^4).
\end{eqnarray*}

The first equation is an irreducible Thue equation, hence we can apply the 
following result by Schmidt:

\vspace{.3cm}

\noindent {\bf Theorem (Schmidt \cite{Schmidt}).--} 
Let $F(x,y)$ be an irreducible binary form of degree $r>3$, with integral 
coefficients. Suppose that not more than $s+1$ coefficients are nonzero. Then the number of solutions of the inequality $|F(x,y)| \leq h$ is, a most, 
$$
(rs)^{1/2} h^{2/r} \left( 1+\log^{1/r}(h) \right).
$$

\vspace{.3cm}

As for our interests are concerned, this gives a bound for the number of possible $(p,q)$ such that
$$
\left| -27(q^4-12q^3p+14q^2p^2+12p^3q+p^4) \right| \leq M.
$$
Hence, as every such solution determines at most one pair in $\cT_5(M)$, 
$$
\left| \cT_5(M) \right| \leq 4 \sqrt{M} \left( 1+\log^{1/4}(M) \right).
$$

A similar result can be applied for points of order $7$. The
equations which must have a solution are now
\begin{eqnarray*}
A &=& -27k^4(p^2-pq+q^2)(q^6+5q^5p-10q^4p^2-15q^3p^3+\\
& & \qquad\qquad\qquad\qquad\qquad\qquad\qquad\qquad\qquad\qquad 30q^2p^4 -11qp^5+p^6), \\
B &=& 54k^6(p^{12}-18p^{11}q+117p^{10}q^2-354p^9q^3+570p^8q^4-486p^7q^5\\
     &   & \quad +273p^6q^6-222p^5q^7+174p^4q^8-46p^3q^9-15p^2q^{10}+6pq^{11}+q^{12}).
\end{eqnarray*}
either for $k=1$ or for $k=1/3$. Hence, using the polynomial defining $B$ and with a similar argument as above
$$
\left| \cT_7(M) \right| \leq 24 \sqrt[6]{M} \left( 1 + \log^{1/12} (M) \right).
$$

Therefore, for all $p$ there is an absolut constant $c_p \in \ZZ_+$ such that
$$
\lim_{M \to \infty} \frac{\left| \cT_p(M) \right|}{|C(M)| }  \leq \lim_{M \to \infty} \frac{c_p M\log(M)}{4M^2-1}
= 0.
$$

This proves the theorem.

\vspace{.3cm}

\noindent {\bf Remark.--} It must be noted here that our arguments are counting pairs $(A,B)$. So, 
in fact, isomorphic curves may appear as separated cases. Both Duke and Grant estimated isomorphism
classes (over $\QQ$) rather than curves. 

But this can also be achieved by the arguments above with a little extra work. We will show now that 
these instances of isomorphic curves are actually negligible as for counting is concerned.

First note that if two curves $E_{(A,B)}$ and $E_{(A',B')}$ are isomorphic over $\QQ$, there must be some 
$u \in \QQ$ such that $A=u^4A'$ and $B=u^6B'$. Hence, there exists some prime $l$ such that, 
say, $l^4|A$ and $l^6|B$ (the case $l^4|A'$ and $l^6|B'$ is analogous). Let us write, for a fixed prime $l$
$$ 
P_n(M,l) = \left\{ x \in \ZZ_+ \; | \; 1\leq x \leq M, \;\; l^n|M \right\},
$$
and by $P_n(M)$ the union of $P_n(M,l)$, where $l$ run the set of prime divisors of $M$.

Then it is clear that
$$
\begin{array}{l}
\displaystyle \left| P_n(M^n) \right| \leq  \displaystyle \sum_{l\leq M} \left| P_n(M^n,l) \right|
= \sum_{l\leq M} \left[ \frac{M^n}{l^n} \right] = \displaystyle \sum_{l\leq M} \left( \frac{M^n}{l^n} + O(1) \right)=\\ 
\quad \,\, =  M^n \sum_{l\leq M} \left( \frac{1}{l^n}\right)  + O(M) = \displaystyle M^n \sum_{l \; \mbox{\scriptsize prime}} \frac{1}{l^n}  + O(M) 
= M^n \mathcal P(n) + O(M),
\end{array}
$$
where $\mathcal P$ is the prime zeta function (see \cite{Froberg}, for instance). So, changing $M^n$ for $M$ we get
\begin{eqnarray*}
\left| P_4(M) \right| &\leq& \displaystyle P(4)M + O \left( \sqrt[4]{M} \right) \simeq 0.0769931M + O \left( \sqrt[4]{M} \right), \\
\left| P_6(M) \right| &\leq& \displaystyle P(6)M + O \left( \sqrt[6]{M} \right) \simeq 0.0170701M + O \left( \sqrt[6]{M} \right). \\
\end{eqnarray*}

Hence, if we are interested in curves up to $\QQ$--isomorphism, our bounds for $|\cT_p(M)|$ are still correct, while we should change 
$$
|\cC(M)| \geq 4M^2-1
$$
by
$$
|\cC(M)| \geq  \left( 4 - P(4)P(6) \right) M^2 + O\left( \sqrt[6]{M} \right)
$$
which obviously makes no difference in the result.

\

\begin{remark} {\rm While all of our boundings for $|\cT_p(M)|$ are of the form $c_p M \log(M)$, computational data show that the actual number of curves on $\cT_p(M)$ depends heavily on $p$, as one might predict after the estimation given by Grant \cite{Grant} for $\cE_p(M)$, the set of elliptic curves $E\in \mathcal{C}_H(M)$ such that $p$ is an excepcional prime for $E$. In fact, a hands--on \verb+Magma+ program gave us the following output
$$
\begin{array}{l|rrrr}
M &  \quad |\cT_2(M)| &   \quad |\cT_3(M)| & 
 \quad  |\cT_5(M)| &  \quad \quad  |\cT_7(M)| \\
\hline
10^4 \,\,\, &  204,220 & 507 & 1 & 1 \\
10^5 &  2,484,196 & 1,935 & 3 & 1 \\
10^6 & 29,430,050 & 5,873 & 11 & 4 \\
10^7 & \,\,\,340,334,782 & 18,387 & 24 & 5 \\
\end{array}
$$

These actual figures are quite smaller than the bounds obtained.
}\end{remark}


\subsection*{Acknowledgement}The first author was supported in part by grants MTM 2009-07291 (Ministerio de Educaci{\'o}n y Ciencia, Spain) and  CCG08-UAM/ESP-3906 (Universidad Auton{\'o}ma de Madrid-Comunidad de Madrid, Spain). The second author was supported by grants FQM--218 and P08--FQM--03894 (Junta de Andaluc\'{\i}a) and MTM 2007--66929 (Ministerio de Educaci\'on y Ciencia, Spain).


\begin{thebibliography}{00}

\bibitem{BI}
M.A. Bennett; P. Ingram: Torsion subgroups of elliptic curves in
short Weierstrass form. Trans. Amer. Math. Soc.  {\bf 357}  (2005)
3325--3337.

\bibitem{Cremona}
J. E. Cremona: Algorithms for modular elliptic curves. Cambridge
University Press, 1992.

\bibitem{Doud}
D. Doud: A procedure to calculate torsion of elliptic curves over
$\QQ$. Manuscripta Math. {\bf 95} (1998) 463--469.

\bibitem{Duke}
W. Duke:  Elliptic curves with no exceptional primes. C.R. 
Acad. Sci. Paris S\'erie I {\bf 325} (1997) 813--818.

\bibitem{Froberg}
C.-E. Fr\"oberg: On the prime zeta function. BIT {\bf 8} (1968) 187--202.

\bibitem{GOT}
I. Garc\'{\i}a--Selfa; M.A. Olalla; J.M. Tornero: Computing the
rational torsion of an elliptic curve using Tate normal form. J.
Number Theory {\bf 96} (2002) 76--88.

\bibitem{GT}
I. Garc\'{\i}a--Selfa; J.M. Tornero: A complete diophantine
characterization of the rational torsion of an elliptic curve.
Available at the arXiv as math.NT/0703578.

\bibitem{GT2}
I. Garc\'{\i}a--Selfa; J.M. Tornero: Thue equations and torsion
groups of elliptic curves. J. Number Theory {\bf 129} (2009) 367--380.

\bibitem{GGT}
I. Garc\'{\i}a--Selfa; E. Gonz\'alez--Jim\'enez, J.M. Tornero:
Galois theory, discriminants and torsion subgroup of elliptic
curves. To appear in Journal of Pure and Applied Algebra.

\bibitem{Grant}
D. Grant: A formula for the number of elliptic curves with exceptional primes. Compositio Math. {\bf 122} (2000) 151--164.

\bibitem{HW}
G.H. Hardy; E.M. Wright: An introduction to the Theory of Numbers 
(5th ed.). Oxford University Press, 1979. 

\bibitem{Husemoller}
D. Husemoller: Elliptic Curves. Springer-Verlag, New York, 1987.

\bibitem{Ingram}
P. Ingram: Diophantine analysis and torsion on elliptic curves.
Proc. London Math. Soc. {\bf 94} (2007) 137--154.

\bibitem{Kubert}
D.S. Kubert: Universal bounds on the torsion of elliptic curves.
Proc. London Math. Soc. {\bf 33} (2) (1976) 193--237.

\bibitem{Lutz}
    E. Lutz: Sur l'equation $y^2=x^3+Ax+B$ dans les corps $p$-adic. J. Reine Angew. Math. {\bf 177} (1937), 431-466.

\bibitem{MazurIHES}
B. Mazur: Modular curves and the Eisenstein ideal. Inst. Hautes
\'Etudes Sci. Publ. Math. {\bf 47} (1977) 33--186.

\bibitem{Mazur}
B. Mazur: Rational isogenies of prime degree. Invent. Math. {\bf
44} (1978) 129--162.

\bibitem{Nagell} 
T. Nagell, Solution de quelque probl\'emes dans la th\'eorie arithm\'etique des cubiques planes du premier genre, Wid. Akad. Skrifter Oslo I, 1935, Nr. 1.


\bibitem{Schmidt}
W.M. Schmidt: Thue equations with few coefficients. Trans. Amer. Math. Soc. {\bf 303} (1987) 241--255.

\bibitem{Serre2}
Serre, J.-P.: Propri\'et\'es galoisiennes des points d'ordre fini des courbes elliptiques. Invent. Math. {\bf 15} (1972), 123--201 (= Collected Papers, III, 1--73).

\bibitem{Silverman}
J.H. Silverman: The arithmetic of elliptic curves. Springer, 1986.


\end{thebibliography}
\end{document}